\documentclass[10pt]{article}

\usepackage{amssymb,amsfonts,graphicx,here,color}
\usepackage[round]{natbib}

\title{\vspace*{1cm} \bf Wavelet analysis of the multivariate fractional Brownian motion}

\author{Jean-Fran\c{c}ois Coeurjolly${}^{1,2}$, Pierre-Olivier Amblard${}^{2,3}$ and Sophie Achard${}^2$  \\
${}^1$ Laboratory Jean Kuntzmann, Grenoble University, France, \\
 ${}^2$ GIPSAlab/CNRS, Grenoble University, France\\
${}^3$ Dept. of  Mathematics\&Statistics, University of Melbourne, Australia.\\
${ }$\\
{\it Corresponding author} \texttt{Jean-Francois.Coeurjolly@upmf-grenoble.fr}}

\newcommand{\R}{\mathbb{R}}
\newcommand{\sign}{\mbox{sign}}

\newtheorem{theorem}{Theorem}
\newtheorem{lemma}[theorem]{Lemma}

\newtheorem{proposition}[theorem]{Proposition}

\newtheorem{theo}{Theorem}

\newtheorem{prop}[theo]{Proposition}

\newcommand\E{\mathbb{E}}

\begin{document}
\def\eqloi{\stackrel{\rm fidi}{=}}

\maketitle

\newpage

\begin{abstract}
The work developed in the paper concerns the multivariate fractional Brownian motion (mfBm) viewed through the lens of the wavelet transform. After recalling some basic properties on the mfBm, we calculate the correlation structure of its wavelet transform. We particularly study the asymptotic behavior of the correlation, showing that if the analyzing wavelet has a sufficient number of null first order moments, the decomposition eliminates any possible long-range (inter)dependence. The cross-spectral density is also considered in a second part. Its existence is proved and its evaluation is performed using a   von Bahr-Essen like representation of the function $\sign(t) |t|^\alpha$. 
The behavior of the cross-spectral density of the wavelet field at the zero frequency is also developed and confirms the results provided by the asymptotic analysis of the correlation. 
\end{abstract}

\centerline{{\bf Keywords}: multivariate fractional Brownian motion, wavelet analysis, cross-correlation, cross-spectrum.}




\section{Motivations and overlook}

The fractional Brownian (fBm) motion developed by \cite{MandVN68} has been extensively studied  as the archetypal model of fractal signals.  Many extensions have been proposed, trying to keep the simplicity of its definition while modeling more complex phenomena. For example, time-dependent Hurst exponent or $d$-dimensional extensions have been introduced which have respectively led to the multifractional Brownian motion \citep{PeltL95} and the fractional Brownian sheet \citep{AyacLP02}. Another extension consists in defining multivariate fractal processes. This extension is needed by many applications ranging from economy to physics, passing by biology and neuroscience \citep{GilA03,AriaC09,AchaBMB08,AchaSWSB06}. In all these disciplines, many modern sensing approaches allow to measure instantaneously different variables from complex phenomena. The need of multivariate signals  models is crucial in order to model and understand these phenomena.

A multivariate extension of the fBm has been proposed recently in a very general setting by \cite{DidiP11} with the help of operator self-similarity. The operator fractional Brownian motion  is an operator self-similar  Gaussian process  with stationary increments. When the operator is diagonal, it is called the multivariate fractional Brownian motion (mfBm). We have particularly studied this diagonal case \citep{AmblCLP11}, 
elaborating on the work of \citet{LavaPS09}. In these works, the correlation structure of the mfBm has been studied. The increments process has also been studied, and we have exhibited its correlation and spectral structures, showing the possible existence of long-range dependence in correlation as well as in cross-correlation between components. 

In this paper, we pursue the study by wavelet analyzing the multivariate process. 
It is now well-accepted that wavelet analysis is the right framework to deal with monovariate fractal signals with stationary increments \citep{Flan88,Flan92,TewfK92,Worn90,FayMRT09,BardLMS00}. Wavelet transforms or decomposition provide a regularized differentiation of the processes, have a filter bank structure in perfect adequacy with $1/f$ type of spectral behavior, may eliminate long-range dependence properties  if the analyzing wavelet is properly chosen. Thus, studying multivariate fractal signals through the lens of the wavelet zoom is indeed natural, and we expect that it will be useful as well in revealing the interaction structure between the components of the mfBm. 

 We thus concentrate on the correlation and on  the spectral structure of the wavelet transform of the mfBm. The principal result of the paper is the explicit form of the second order statistics of this multivariate Gaussian random field.  We study the asymptotic behavior of the cross-correlation function of the wavelet, and this allows us to prove that choosing a wavelet with at least two null first order moments eliminates any possible long dependence in the correlations. We  prove the existence and calculate the cross-spectral density of the wavelet field. We thus extend the result of \cite{KatoM99} providing the existence of the spectral density of the wavelet transform of the fBm. The proof uses a generalization of the von Bahr-Essen representation of $|t|^\alpha$ (used by \citet{KatoM99}) to the function $\sign(t) |t|^\alpha$. We also provide the behavior of the density at the zero frequency, corroborating the asymptotic result obtained in the time domain. We stress on the fact that our asymptotic results do not impose that the wavelet is a real function, nor that it has a compact support. 

The paper is organized as follows. In order to have a self-contained exposition, we recall some basics   definition and results on the mfBm in Section~1. In Section~2, we set the wavelet analysis and look at some self-similarity properties inherited from the process. The full correlation structure is developed in Section~3 whereas the spectral counterpart is in Section~4. The needed generalized von Bahr-Essen representation is proved in the last section.

\section{Some facts on the multivariate fractional Brownian motion}

The $p$ dimensional multivariate fractional Brownian motion (mfBm) $x(t)$ is defined as the only Gaussian process having stationary increments and having components jointly self-similar with parameters $(H_1,\ldots,H_p) \in (0,1)^p$. The self-similarity property can be stated as follows: for any real $\lambda > 0$, $x(\lambda t) \eqloi \lambda^H x(t)$ where $H=\mbox{diag}(H_1,\ldots,H_p) $ and $\lambda^H$ is intended in the matrix sense. The notation $\eqloi$ stands for equality of all the finite-dimensional probability distributions. 

The cross-covariance structure induced by the multivariate self-similarity property and the stationarity of the increments has been first studied by \citet{LavaPS09}, Theorem~2.1, without having recourse to the Gaussian assumption. \citet{AmblCLP11} have parameterized this covariance structure in a more simple way as follows.

\begin{prop}[Proposition~3, \citep{AmblCLP11}] Let $j,k\in \{1,\ldots,p\}$, $j\neq k$, then there exists $\sigma_j>0$, $(\rho_{jk},\eta_{jk})\in [-1,1] \times \R$ satisfying $\rho_{jk}=\rho_{kj}$  and $\eta_{jk}=-\eta_{kj}$, such that
\begin{equation}  \label{covariancedef:eq}
r_{jk}(s,t) := \E[  x_j(s) x_k(t) ] =
  \frac{\sigma_j\sigma_k}{2} \left\{ w_{jk}(-s) +
  w_{jk}(t)    - w_{jk}(t-s) \right\},
\end{equation}
where the function $w_{jk}(h)$ is defined  by
\begin{equation}
w_{jk}(h) = \left\{ \begin{array}{ll}
(\rho_{jk} -\eta_{jk} \sign (h) ) |h|^{H_j+H_k} & \mbox{ if } H_j+H_k \neq 1, \\
{\rho}_{jk} |h| + {\eta}_{jk} h \log|h| & \mbox{ if } H_j+H_k =1.
\end{array} \right.
\label{fonctionw:eq}
\end{equation}
\end{prop}
This result is also valid in the case $j=k$ when setting $\rho_{jj}=1$ and $\eta_{jj}=0$ in~(\ref{fonctionw:eq}); we thus recover  the covariance structure of a monovariate fBm. The parameter $\sigma_j^2$ is the variance of a fBm at time $1$,  $\mbox{Var}(x_j(1))$, whereas $\rho_{jk}$ represents the instantaneous correlation between components $j$ and $k$ at time 1, i.e. $\E [x_j(1)x_k(1)]$. The antisymmetric parameter $\eta_{jk}$ is related to the time-reversibility property of the multivariate process. Indeed, the mfBm is time reversible, i.e. $x(t) = x(-t)$ in distribution for every $t$, if and only if $\eta_{jk}=0$ for all $j,k$ \citep{AmblCLP11}. 

\begin{center}
\begin{figure}[h]
\begin{tabular}{lll}
\includegraphics[width=5cm]{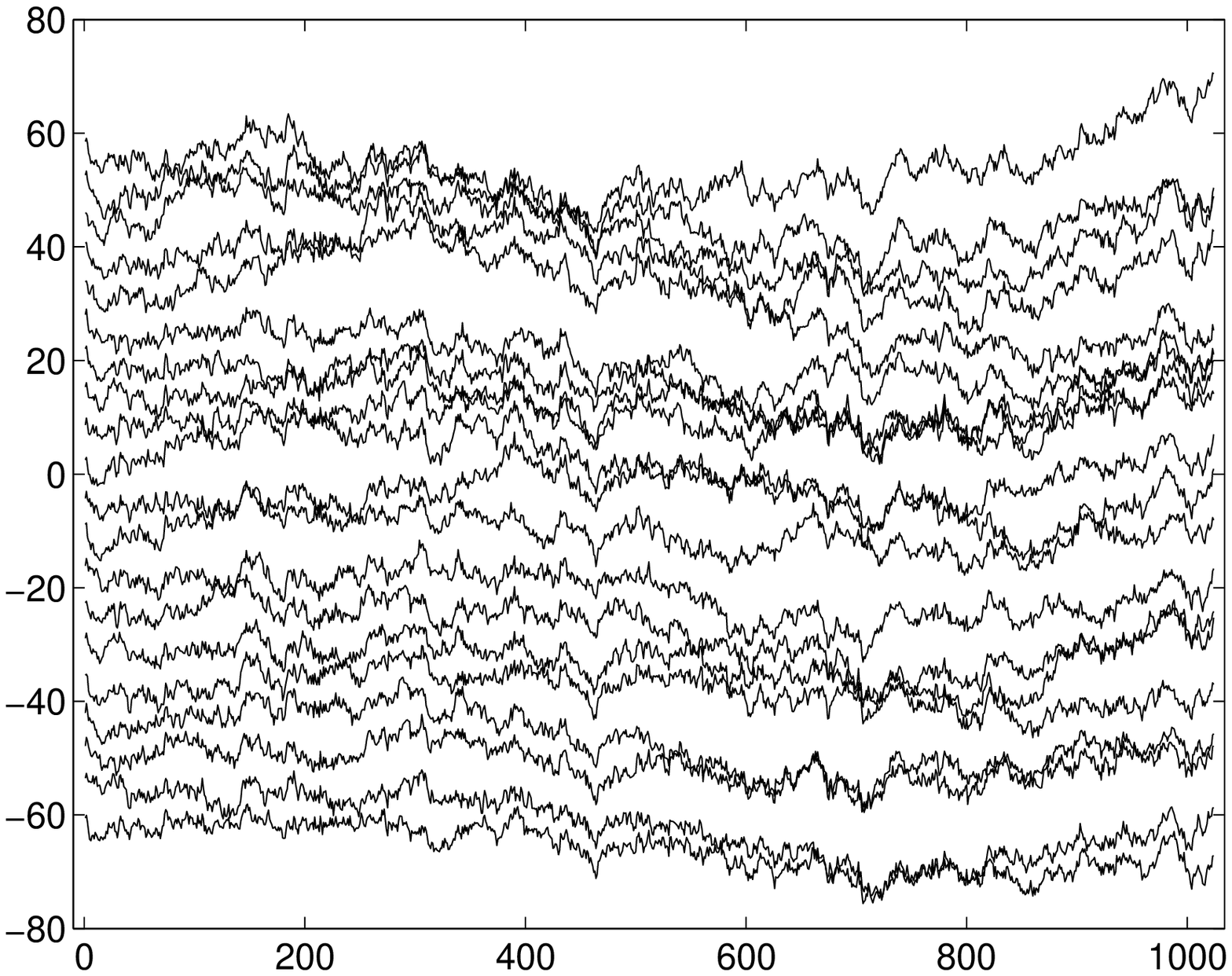} &\includegraphics[width=5cm]{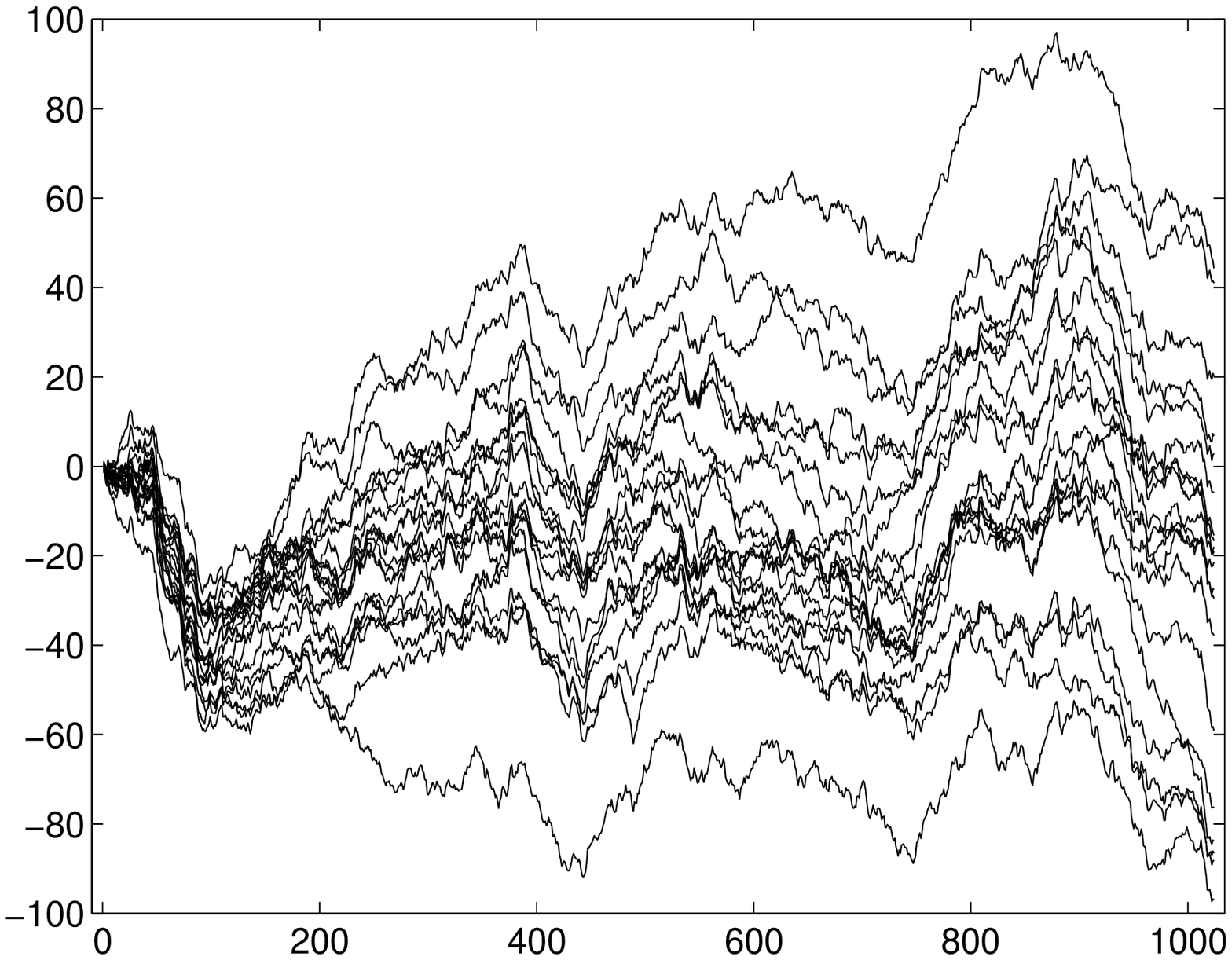}&\includegraphics[width=5cm]{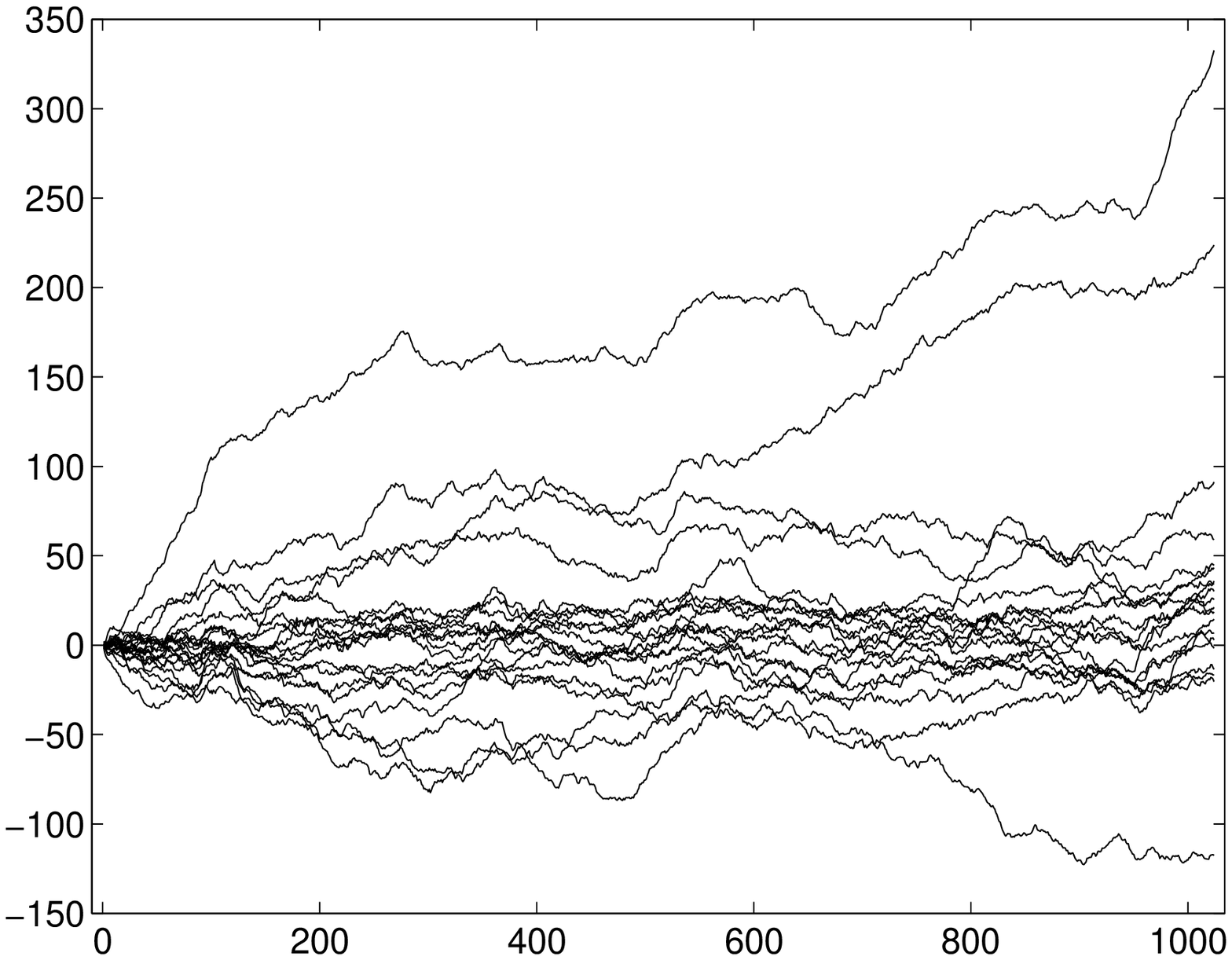}
\end{tabular}
\caption{\small Examples of discretized sample paths of a time reversible ($\eta_{jk}=0$) mfBm of length $n=1024$, with $p=20$ components. The Hurst exponents are equally spaced in $[0.3 , 0.4]$ (left plot), $[0.6 , 0.7]$ (middle plot) and $[0.4, 0.8]$ (right plot). The correlation parameters are set to 0.7 (left and middle plot) and to 0.3 (right plot). The components are shifted artificially in the left plot for the sake of visibility.} 
\label{ex_mfbm:fig}
\end{figure}
\end{center}

To ensure that the matrix  given by (\ref{covariancedef:eq}) is the cross-covariance matrix of a process, the constraints imposed on $\rho_{jk}$ and $\eta_{jk}$ are not sufficient. A necessary and sufficient condition, proved by \citet{AmblCLP11}, corresponds to the positive-definiteness of the Hermitian matrix with entries $\Gamma(H_j+H_k+1) \times \xi_{jk}$ where $\xi_{jk}$ is defined by
\begin{equation} \label{condExist}
\xi_{jk}= \left\{ \begin{array}{ll}
\rho_{jk}\sin\left( \frac\pi2(H_j+H_k) \right) - \mathbf{i} \eta_{jk} \cos\left( \frac\pi2(H_j+H_k) \right) & \mbox{ if } H_j+H_k \neq 1\\
\rho_{jk} - {\bf i} \frac\pi2 \eta_{jk} & \mbox{ if } H_j+H_k =1,
\end{array} \right.
\end{equation}
where ${\bf i}=\sqrt{-1}$. For example, when $p=2$ and $\eta_{jk}=0$ no condition is required for the correlation $\rho_{12}$ but when $H_1=H_2$, and when $H_1=0.1$ and $H_2=0.2$ the correlation $\rho_{12}$ cannot exceed $0.514$, see \citep{AmblCLP11} for more discussion. The problem of simulation of such a process has been investigated in \citet{AmblCLP11} using the   algorithm developed by \citet{ChanW99}. Figure (\ref{ex_mfbm:fig}) presents some examples in order to illustrate the process.

The covariance structure of the   increment process (at lag 1) can be easily deduced from~(\ref{covariancedef:eq}). When $j=k$, we obviously recover the covariance of the fractional Gaussian noise and the classical property that this process is short-memory when $H_j\leq 1/2$ and long-memory when $H_j>1/2$. In the multivariate case, long-range (interdependence) may also appear in the cross-covariance. Indeed~\cite{AmblCLP11} proved that for all $j\neq k$, the cross-covariance behaves asymptotically as $|h|^{H_j+H_k-2}$ (up to a constant) meaning that the long-memory property arises as soon as $H_j+H_k \geq 1$ which can appear in three different situations: $H_j=H_k=1/2$, $H_j<1/2$ and $H_k>1-H_j$ or $H_j>1/2$ and $H_k>1/2$. In those cases, some troubles may appear when it comes to infer parameters of the models from data. Indeed, long-range dependence may lead to very slow convergence of estimators. As already observed in many works, {\it e.g.} \citep{Flan92,VeitA99,BardLMS00,Coeu01}, having recourse  to wavelet types of transformation is an elegant way to overcome the problem. Indeed, using wavelet types of transformation with a correctly chosen filter allows to extract the stationary part from the fBm and allows  to ``whiten" the increments. We describe such an approach in the following section.

\section{Wavelet Analysis: definition, stationarity and self-similarity}
\label{wavelet:sec}

The use of wavelet analysis in the understanding of the monovariate fractional Brownian motion, and more generally for  the study of fractal signals, goes back to the early works of  \citet{Flan88,Flan92},  \citet{TewfK92},  \citet{Worn90} to cite some but a few. 

The aim is now to analyze the multivariate fractional Brownian motion through the lens of the wavelet transform. We use the continuous wavelet transform here, but a similar analysis could be performed in  the multiresolution framework using orthonormal wavelet bases.  We will consider complex valued wavelets, not necessarily in the Hardy class, not necessarily with compact support. The hypothesis we impose on the wavelets will be detailed when needed. 

\subsection{Definition and stationarity}

Let $\psi$ be a  complex wavelet function, let $a>0$ and $b\in \R$ and consider $\psi_{ab}(.) =a^{-1/2} \psi((.-b)/a)$. Let
\begin{eqnarray}
d^j_{a,b} &:=& \left< x_j \Big| \psi_{a,b} \right>_{L^2} = \int_{\R} x_j(t) \overline{\psi_{a,b}}(t) dt \label{djab}
\end{eqnarray}
the wavelet transform of the $j$th component of a multivariate fractional Brownian motion.  $\overline{\psi}$ denotes the complex conjugate of $\psi$. In this section, we assume that conditions [C1] and [C2(2)] are satisfied, where:\\

\noindent{[C1]} Admissibility condition: $\psi(t) \in L^2$ and $ |\widehat{\psi}(\omega)|^2/|\omega| \in L^1$, where $\widehat{\psi}$ is the Fourier transform of~$\psi$. \\

\noindent {[C2(K)]} $t^m  \psi(t) \in L^1$ for $m=0,1,\ldots,K$.\\

Condition [C1] ensures that $\widehat{\psi}(0)=0$ and that $\int_\R \psi(t)dt=0$. We note, as \citet{KatoM99}, that under condition [C2(1)], the integral (\ref{djab}) is well-defined as a sample path integral and is a second-order random variable. This follows, since under [C2(1)] we have $\int_{\R} |s|^H |\psi_{a,b}(s)| ds<+\infty, \forall H \in (0,1)$.

 The aim of this section is to focus on the correlation between the wavelet transforms (at different scales and different times) of two components $j$ and $k$ of the multivariate fractional Brownian motion. The wavelet transform is a random field. It is clearly zero mean and Gaussian. We have for $a_1,a_2>0$ and $b,h \in \R$
\begin{eqnarray*}
E[d^j_{a_1,b+h}\overline{d^k_{a_2,b}} ]  = \int_{\R^2} r_{jk}(t_1,t_2) \overline{\psi_{a_1,b+h}}(t_1) \psi_{a_2,b}(t_2)  dt_1 dt_2.
\end{eqnarray*}
Under [C1], and from (\ref{fonctionw:eq}) the last expression reduces to 
\begin{eqnarray*}
E[d^j_{a_1,b+h}\overline{d^k_{a_2,b}} ]  = -\frac{\sigma_j \sigma_k }2 \int_{\R^2}  w_{jk}(t_2-t_1)\overline{\psi_{a_1,b+h}}(t_1) \psi_{a_2,b}(t_2)  dt_1 dt_2.
\end{eqnarray*}
Let $\Gamma_\psi(v):=\int_\R \psi_{a_1,b+h}(u) \overline{\psi_{a_2,b}}(u+v) du$ be the correlation function between the two wavelets $\psi_{a_1,b+h}$ and $\psi_{a_2,b}$. Then we have
\begin{eqnarray}
E[d^j_{a_1,b+h}\overline{d^k_{a_2,b}} ]  = - \frac{\sigma_j \sigma_k }2\int_\R w_{jk}(v) \overline{ \Gamma_\psi}( v)  dv. \label{eq-wjkgpsi}
\end{eqnarray}
Note that [C2(2)] implies that for all the values of $H_j$ and $H_k$,  $\int_\R |w_{jk}(v)| \left| \Gamma_\psi(v)\right| dv<+\infty$. With two changes of variables, this may also be rewritten as
\begin{eqnarray}
E[d^j_{a_1,b+h}\overline{d^k_{a_2,b}} ] &=&  - \frac{\sigma_j \sigma_k }2\ \sqrt{a_1a_2}\times
\int_{\R^2}\!\! w_{jk}(a_2t_2-a_1t_1-h) \overline\psi(t_1)\psi(t_2)dt_1dt_2. \label{cor-wt}
\end{eqnarray}
If we  interpret for fixed parameters $a_1$ and $a_2$, the quantity $E[d^j_{a_2,b+h}\overline{d^k_{a_2,b}} ]$ as the cross-correlation between two signals, we observe that it depends only on the difference between the times at which it is evaluated (i.e. $h$). With the fact that the wavelet transform is a zero mean and Gaussian field, we conclude that $d^j_{a_1,.}$ and $d^k_{a_2,.}$ are jointly stationary signals. 

\subsection{Self-similarity type property of the cross-wavelet transform}

The variance of the wavelet transform at similar scales for the fractional Brownian motion with Hurst parameter $H$ exhibits some self-similarity. Indeed, it is proved in \citet{Flan88} for example that for all $b$
$$
Var(d^j_{a,b} ) = a^{2H+1} \times \left(-\frac{\sigma^2}2\int_{\R^2} |t_2-t_1|^{2H}  \overline\psi(t_1)\psi(t_2)dt_1dt_2\right).
$$
We note here that the same behavior holds for the cross-wavelet variance. 

\begin{proposition} Under the assumptions [C1] and [C2(2)], let $h=0$ and fix $a_1=a_2=a>0$. Then, 
\begin{equation} \label{ss-wt}
E[d^j_{a,b}\overline{d^k_{a,b}} ] = a^{H_j+H_k+1}\; \left(- \frac{\sigma_j \sigma_k }2 \; z_{jk} \right) 
\quad \mbox{ and } \quad 
 Corr[d^j_{a,b}, {d^k_{a,b}} ] =\frac{z_{jk}}{\sqrt{z_{jj} z_{kk} }},\nonumber
\end{equation}
where $z_{jk}:= \int_{\R^2} w_{jk}(t_2-t_1) \overline\psi(t_1)\psi(t_2)dt_1dt_2$.
\end{proposition}

\begin{proof}
Consider Equation~(\ref{cor-wt}). The result is obvious when $H_j+H_k\neq1$ since for any $a>0$, $w_{jk}(av)=a^{H_j+H_k} w_{jk}(v)$. Now, when $H_j+H_k=1$, the result comes from Condition [C1] ensuring that $\int_{\R^2}\eta_{jk} \times a(t_2-t_1)\log(a) \overline\psi(t_1)\psi(t_2)dt_1dt_2 = 0$.
\end{proof}

Let us observe that the instantaneous cross-wavelet correlation is independent of the scale. 


\section{Cross-correlation structure of the wavelet transform of the mfBm}

For fixed scales, $a_1, a_2$, we now specify the behavior of the cross-wavelet covariance (or correlation) as $|h|\to +\infty$. In particular, our aim is to exhibit the influence of the number of vanishing moments of the wavelet function on the asymptotic cross-wavelet covariance. Such a result needs the following assumption:\\

\noindent {[C3]} The wavelet function has $M\geq 1$ vanishing moments that is 
$$\int_\R t^m \psi(t)dt=0 \mbox{ for } m=0,\ldots,M-1 \qquad \mbox{ and } \qquad \int_\R t^M\psi(t)dt \neq 0.$$
We may now derive our result obtained as $|h|\to +\infty$. Let us first recall Landau notation: for two functions $f(h)$ and $g(h)$ defined on $\R$, we denote by $f(h) {\sim} g(h)$ as $|h|\to +\infty$ (resp. $f(h){=}o(g(h))$ and $f(h){=}\mathcal{O}(g(h))$) if $\lim_{|h|\to +\infty} f(h)/g(h)=1$ (resp. $\lim_{|h|\to +\infty} f(h)/g(h)=0$ and  $f(h)/g(h)$ is bounded for all $h$).

\begin{theorem} \label{prop-corwav} 
Assume [C1], [C2(2M+1)] and [C3] hold, then as $|h|\to +\infty$, we have
$$
E[d^j_{a_1,b+h}\overline{d^k_{a_2,b}} ] \sim -\frac{\sigma_j \sigma_k }2 \kappa(\psi,M) \;  {\tau}_{jk}(h)\;|h|^{H_j+H_k-2M} 
$$
where $\kappa(\psi,M):={2M \choose M} (a_1 a_2)^{M} \left| \int t^M \psi(t)dt \right|^2$ and
\begin{equation}\label{def-taut}
{\tau}_{jk}(h) = \left\{
\begin{array}{ll}
(\rho_{jk} + \eta_{jk}\; \sign(h)) \mbox{${ {H_j+H_k \choose 2M}}$} & \mbox{ if } H_j+H_k\neq 1\\
-\frac{\eta_{jk} \times \sign(h)}{2M(2M-1)} & \mbox{ if } H_j+H_k=1.
\end{array} \right.
\end{equation}
\end{theorem}
We notice that the equivalence stated has a meaning as soon as $\rho_{jk}+\eta_{jk}\; \sign(h)\neq 0$  when $H_j+H_k\neq 1$ and as soon as $\eta_{jk}\neq 0$ when $H_j+H_k=1$. In the opposite cases, a careful look at the proof shows that the equivalence can be replaced by an upper-bound or more precisely $\E[ d^j_{a_1,b+h}\overline{d^k_{a_2,b}} ]= o(|h|^{1-2M})$.\\


\vspace{.5cm}

Similarly to the fractional Brownian motion, Theorem~\ref{prop-corwav} asserts that the higher $M$,  the less correlated the wavelet transforms of the components $j$ and $k$ of the multivariate fractional Brownian motion. This has many implications. In particular, this suggests that estimating the instantaneous cross-wavelet correlation at a scale $a$ may be efficiently done by using the empirical correlation since at scale $a$, $d_{a,b+h}^j$ and $d_{a,b}^k$ are not too much correlated if $M$ is large.\\

\begin{proof}
The proof is split into two cases. Before this, we denote by ${D_h}:=\{ (t_1,t_2)\in \R^2: |a_2t_2-a_1t_1|< \frac{|h|}{2} \}$ and we note in particular that 
$$
\begin{array}{ll}
\forall (t_1,t_2) \in D_h & \left| \frac{a_2 t_2 -a_1t_1}h \right| \leq \frac12 <1 \quad \mbox{ and } \quad
\sign ( a_2t_2-a_1 t_1 - h ) = -\sign(h), \\
&\\
\forall (t_1,t_2) \in \R^2\setminus D_h & \left| \frac{h}{a_2 t_2 -a_1t_1} \right| \leq 2.
\end{array}
$$

\noindent \underline{Case 1}. $\alpha:=H_j+H_k\neq 1$. \\

We assume here that $\rho_{jk} + \eta_{jk} \sign(h) \neq 0$. Let us write
$
E[d^j_{a,b+h}\overline{d^k_{a,b}} ] = -\frac{\sigma_j \sigma_k }2 \sqrt{a_1 a_2}\times T
$
with
$$
T:=\int_{\R^2} (\rho_{jk}-\eta_{jk} \;\sign (a_2t_2-a_1 t_1 - h) |a_2t_2-a_1 t_1 - h|^\alpha \overline{\psi}(t_1) \psi(t_2) dt_1 dt_2 =T_1+T_2, 
$$
and where $T_1$ (resp. $T_2$) corresponds to the integral on ${D_h}$ (resp. ${\R^2\setminus {D_h}}$). Let us first prove that 
$|h|^{2M- \alpha} T_2 \to 0$ as $|h|\to +\infty$. Denoting $c^\vee = |\rho_{jk}|+|\eta_{jk}|$, we have (since $2M-\alpha>0$)
\begin{eqnarray*}
|h|^{2M- \alpha} |T_2| &\leq & 
c^\vee\int_{\R^2\setminus {D_h}} \!\!\!|a_2t_2-a_1t_1|^{\alpha}|h|^{2M-\alpha}  \left| 1+ \frac{h}{a_2t_2-a_1t_1}\right|^{\alpha} |\psi(t_1)| |\psi(t_2)| dt_1dt_2 \\
&\leq &2^{2M-\alpha}  3^{\alpha} c^\vee \int_{\R^2\setminus {D_h}} \!\!\!(a_2t_2-a_1t_1)^{2M} |\psi(t_1)| |\psi(t_2)| dt_1dt_2.
\end{eqnarray*}
The result is then obtained by using assumption [C2(2M)] and the dominated convergence theorem. Now, within the domain ${D_h}$, one may use the series expansion of $(1+x)^\alpha$ (for $|x|<1$).
\begin{eqnarray*}
T_1 &=& |h|^\alpha \int_{D_h} (\rho_{jk} -\eta_{jk} \; \sign(a_2t_2-a_1t_1-h)) \left(1- \frac{a_2t_2-a_1t_1}h \right)^\alpha \overline{\psi}(t_1) \psi(t_2) dt_1 dt_2 \\
&=& |h|^\alpha (\rho_{jk} +\eta_{jk} \; \sign(h)) \int_{D_h} \left( \sum_{\ell \geq 0} (-1)^\ell{\alpha \choose \ell}\left(\frac{a_2t_2-a_1t_1}{h}\right)^\ell \right) \overline{\psi}(t_1) \psi(t_2) dt_1 dt_2,
\end{eqnarray*}
where ${\alpha \choose \ell}$ denotes the binomial coefficient $(\alpha)(\alpha-1)\ldots(\alpha-\ell+1)/\ell!$. Decompose $T_1$ into three terms (denoted by $T_1^\prime,T_2^\prime$ and $T_3^\prime$) corresponding to the $2M$ first terms of the series, the $(2M+1)$th term ($\ell=2M$) and the remainder terms. Then,
$$
T_1^\prime = |h|^{\alpha} (\rho_{jk}+\eta_{jk}\;\sign(h))\sum_{\ell=0}^{2M-1}(-1)^\ell h^{-\ell}{\alpha \choose \ell} \int_{D_h} \left(a_2t_2-a_1t_1\right)^\ell \overline{\psi}(t_1) \psi(t_2) dt_1 dt_2. 
$$
Under Assumption [C3], $\psi$ has $M$ vanishing moments and therefore the previous expression reduces to
$$
T_1^\prime = - |h|^\alpha (\rho_{jk}+\eta_{jk}\;\sign(h))\sum_{\ell=0}^{2M-1}(-1)^\ell h^{-\ell}{\alpha \choose \ell} \int_{\R^2\setminus {D_h}} \left(a_2t_2-a_1t_1\right)^\ell \overline{\psi}(t_1) \psi(t_2) dt_1 dt_2. 
$$
Now,
\begin{eqnarray*}
|h|^{2M-\alpha} |T^\prime_1| &\leq&  c^\vee\sum_{\ell=0}^{2M-1}\left|{\alpha \choose \ell}\right| \int_{\R^2\setminus {D_h}} 2^{2M-\ell}\left(a_2t_2-a_1t_1\right)^{2M} |{\psi}(t_1)| |\psi(t_2)| dt_1 dt_2.
\end{eqnarray*}
Assumption [C2(2M)] and the dominated convergence theorem may be combined to prove that $|h|^{2M-\alpha} T^\prime_1\to 0$ as $|h|\to +\infty$. The term $T^\prime_2$ is defined as
$$
T^\prime_2 := |h|^{\alpha-2M} (\rho_{jk}+\eta_{jk}\;\sign(h)){\alpha \choose 2M} \int_{D_h}  \left(a_2t_2-a_1t_1\right)^{2M} \overline{\psi}(t_1) \psi(t_2) dt_1 dt_2. 
$$
As previously we obtain 
\begin{eqnarray*}
\frac{|h|^{2M-\alpha} T^\prime_2}{\rho_{jk}+\eta_{jk}\;\sign(h)}  &\to& {\alpha \choose 2M} \int_{\R^2}  \left(a_2t_2-a_1t_1\right)^{2M} \overline{\psi}(t_1) \psi(t_2) dt_1 dt_2 \\
&=& {\alpha \choose 2M} {2M \choose M} (a_1 a_2)^{M} \left| \int t^M \psi(t)dt \right|^2 =  {\alpha \choose 2M}\kappa(\psi,M).
\end{eqnarray*}
Since $T=T_1+T_2 = T^\prime_1+T^\prime_2+T^\prime_3+T_2$, the proof will be completed if we manage to prove that $|h|^{2M-\alpha}T^\prime_3\to 0$. Let us  write
\begin{eqnarray*}
|h|^{2M-\alpha} T^\prime_3 &=& h^{2M} (\rho_{jk}+\eta_{jk} \; \sign(h))\int_{D_h}   \sum_{\ell\geq 2M+1}  (-1)^\ell {\alpha \choose \ell}  \left( \frac{a_2t_2-a_1t_1}{h} \right)^\ell \overline{\psi}(t_1) \psi(t_2) dt_1 dt_2 \\
&=& \frac{\rho_{jk}+\eta_{jk} \; \sign(h) }{h} \;\times \\
&&\int_{D_h}  (a_2t_2-a_1t_1)^{2M+1} \left(\sum_{\ell\geq 0} (-1)^{\ell+1} {\alpha \choose \ell+2M+1} \left( \frac{a_2t_2-a_1t_1}{h} \right)^\ell \right)\;\;\overline{\psi}(t_1)\psi(t_2) dt_1 dt_2.
\end{eqnarray*}
The binomial coefficient appearing in the last equation satisfies, with $\ell'=\ell+2M+1$
\begin{eqnarray*}
\Big|  {\alpha \choose \ell'} \Big| &=& \frac{\big| \alpha (\alpha-1) \cdots (\alpha-\ell'+1) \big|}{\ell'!} \\
&\leq& \frac{ 2  ( 2- \alpha) \cdots (\ell'-1-\alpha)}{\ell'!}  \quad\mbox{ since } \alpha \leq 2 \\
&\leq & \frac{ 2 (\ell' -1)!}{\ell'!} = \frac{ 2 }{\ell'} \leq \frac2\ell.
\end{eqnarray*}
Recall that in $D_h$ we have $\big|a_2t_2-a_1t_1\big|/|h|\leq 1/2$. The series in the previous integral then satisfies
\begin{eqnarray*}
\Big|\sum_{\ell\geq 0}  {\alpha \choose \ell+2M+1} \left( \frac{a_2t_2-a_1t_1}{h} \right)^\ell\Big| &\leq&\Big| {\alpha \choose 2M+1}\Big| + \sum_{\ell\geq 1}  \left|{\alpha \choose \ell+2M+1}\right| \left| \frac{a_2t_2-a_1t_1}{h} \right|^\ell\\
&\leq& \frac{2}{2M+1} + \sum_{\ell\geq 1} \frac{2}{\ell} \; 2^{-\ell} \\
&=&  \frac{2}{2M+1} +2 \log(2) =: C_M.   
\end{eqnarray*}
Thus we obtain 
\begin{eqnarray*}
|h|^{2M-\alpha} \big|T^\prime_3\big| &\leq& \frac{C_M c^\vee}{|h|} \int_{\R^2}  \big|a_2t_2-a_1t_1\big|^{2M+1}  |\psi(t_1)| |\psi(t_2)| dt_1 dt_2.
\end{eqnarray*}
Since by Assumption [C2(2M+1)], $t^{2M+1} \psi(t)\in L^1$, we have $|h|^{2M-\alpha} \big|T^\prime_3\big| =O(|h|^{-1})$, whence the result. \\

\noindent \underline{Case 2}. $H_j+H_k =1$.

We assume here that $\eta_{jk}\neq0$. We take the same notation as previously. We first note that, under [C1], the term $T$ can be rewritten as
$$
T =\int_{\R^2} \rho_{jk}|a_2t_2-a_1t_1-h|+ \eta_{jk} (a_2t_2-a_1t_1-h) \log\left|1-\frac{a_2t_2-a_1t_1}{h} \right| \overline{\psi}(t_1)\psi(t_2) dt_1 dt_2. 
$$ 
We decompose $T$ in $T_1+T_2$ (as done in case 1). The proof that $|h|^{2M-1} T_2\to 0$ as $|h|\to +\infty$ follows similar arguments as in the case 1 and is therefore omitted. Now, the term $T_1$ can be rewritten as
\begin{eqnarray*}
T_1 &=& \rho_{jk}|h| \int_{D_h} \left(1- \frac{a_2t_2-a_1t_1}{h}\right) \overline{\psi}(t_1)\psi(t_2) dt_1 dt_2 \\
&&- \eta_{jk} h \int_{D_h} \left(1- \frac{a_2t_2-a_1t_1}{h}\right) \log\left(1- \frac{a_2t_2-a_1t_1}{h}\right) \overline{\psi}(t_1)\psi(t_2) dt_1 dt_2.
\end{eqnarray*}
Denote by $\tilde{T}_1$ and  $\tilde{T}_2$ these two terms. Assumption [C1] leads to 
\begin{eqnarray*}
\tilde{T}_1&=&-\rho_{jk} |h| \int_{\R^2\setminus {D_h}} \left(1- \frac{a_2t_2-a_1t_1}{h}\right) \overline{\psi}(t_1)\psi(t_2) dt_1 dt_2\\
&=& -\rho_{jk} \int_{\R^2\setminus D_h} |a_2t_2-a_1t_1| \left| 1 - \frac{h}{a_2t_2-a_1t_1}\right| dt_1dt_2. 
\end{eqnarray*}
Then, we assert that
\begin{eqnarray*}
|h|^{2M-1} |\tilde{T}_1|& \leq & 2^{2M-1} 3 |\rho_{jk}| \int_{\R^2\setminus {D_h}} \Big(a_2t_2-a_1t_1\Big)^{2M} |{\psi}(t_1)| |\psi(t_2)| dt_1 dt_2 \to 0
\end{eqnarray*}
as $|h|\to +\infty$. For the term $\tilde{T}_2$, we may use the series expansion of $\log(1+x)$ (for $|x|<1$). We omit the details and leave the reader to verify that as $|h|\to +\infty$
\begin{eqnarray*}
\tilde{T}_2 &\sim& \eta_{jk} h \int_{\R^2}  \left( 1- \frac{a_2t_2-a_1t_1}{h}\right) \left( \frac{-1}{2M-1} \left(  \frac{a_2t_2-a_1t_1}{h}\right)^{2M-1} \right. \\
&&\left.  \frac{-1}{2M} \left(  \frac{a_2t_2-a_1t_1}{h}\right)^{2M}  \right)\overline{\psi}(t_1)\psi(t_2) dt_1 dt_2 \\
&\sim& -h^{1-2M} \times \frac{\eta_{jk}}{2M(2M-1)} {2M \choose M} (a_1 a_2)^{M} \left| \int t^M \psi(t)dt \right|^2.
\end{eqnarray*}
Hence, $T\sim |h|^{1-2M} \times \left( -\frac{\eta_{jk} \times \sign(h)}{2M(2M-1)} \right)\kappa(\psi,M)$. 

In this proof, Fubini's theorem and interchanges of integrals and (in)finite sums are widely used. All of these are justified by the absolute convergence of the different series related to the expansions of $(1+x)^{\alpha}$ or $\log(1+x)$ for $|x|<1$ and Assumption~[C2(2M+1)].
\end{proof}


\section{Cross-spectral density of the wavelet transform of the mfBm} \label{sec-crosswavspd}

In the case of the fBm, the expression of the spectral density of the wavelet transform was provided by \citet{Flan88,Flan92}. A rigorous proof of the existence of this spectral density in the $L^1$ sense was obtained by \citet{KatoM99}. On the basis of this work, our ambition is to provide the cross-spectral density between wavelet transforms (at different scales) of components $j$ and $k$ of the multivariate fractional Brownian motion.
The idea is to obtain the following spectral representation for the cross-correlation
\begin{eqnarray*}
E[d^j_{a_1,b+h}\overline{d^k_{a_2,b}} ] = \frac{1}{2\pi } \int_\R {S}_{a_1,a_2}^{jk}(\omega) e^{{\bf i} \omega t} d\omega.
\end{eqnarray*}


\begin{theorem} \label{prop-spdwt} Under Assumptions [C1], [C2(M)] and [C3] (with $M\geq 2$), we derive the following assertions.\\
$(i)$ The cross-spectral density of the wavelet transforms of two components $j$ and $k$ 
exists and is given by
\begin{eqnarray} 
{S}_{a_1,a_2}^{jk}(\omega)  &=&\sqrt{a_1a_2}  \; {\sigma_j \sigma_k \Gamma(H_j+H_k+1)}   \;\zeta_{jk}(\omega)\; \frac{\overline{\widehat{\psi}(a_1\omega)} {\widehat{\psi}(a_2\omega)} }{|\omega|^{H_j+H_k+1}}
 \label{eq-spdwt}
\end{eqnarray}
where 
$$ \zeta_{jk}(\omega)= \left\{ \begin{array}{ll}
\rho_{jk} \sin\left(\frac{\pi}2 (H_j+H_k)\right)\; + \;{\bf i}\; \eta_{jk} \cos\left(\frac{\pi}2 (H_j+H_k)\right) \sign(\omega) & \mbox{ if } H_j+H_k \neq 1 \\
\rho_{jk} + {\bf i} \frac\pi2 \eta_{jk} \;\sign(\omega) &  \mbox{ if } H_j+H_k= 1. \\
\end{array} \right.
$$

$(ii)$ We have for both cases, as $\omega\to 0$
$$
\Big| {S}_{a_1,a_2}^{jk}(\omega) \Big| \sim  (a_1a_2)^{M+1/2} \sigma_j\sigma_k\; \Gamma(H_j+H_k+1) \; |\widehat{\psi}^{(M)}(0)|^2 \;|\zeta_{jk}(\omega)| \;      |\omega|^{2M-1-\alpha}.
$$
$(iii)$ Moreover, the coherence function between the two components $j$ and $k$ satisfies:
\begin{equation}\label{eq-coherence:wt}
{C}_{a_1,a_2}^{jk}(\omega) : = \frac{\left| {S}_{a_1,a_2}^{jk}(\omega) \right|^2}{{S}_{a_1,a_1}^{jj}(\omega) {S}_{a_2,a_2}^{kk}(\omega) } = |\zeta_{jk}(\omega)|^2 \; \frac{\Gamma(H_j+H_k+1)^2}{\Gamma(2H_j+1)\Gamma(2H_k+1)} \times
 \frac{{\widehat{\psi}(a_1\omega)} \overline{\widehat\psi(a_2 \omega)}}{\overline{\widehat{\psi}(a_1\omega)} {\widehat\psi(a_2 \omega)}}.
\end{equation}
\end{theorem}

Before writing down the proof, let us give some comments.
\begin{enumerate}
\item Item (ii) in Theorem~\ref{prop-spdwt} is the spectral analogue of Theorem~\ref{prop-corwav}. Indeed the behavior ot the cross-correlation at infinite lags is linked to the behavior of its Fourier transform at the zero frequency.  We recover the fact that as soon as $M> H_j+H_k+1/2$, the long-range interdependence is destroyed. The divergence of  $ |\omega|^{-1-H_j+H_k} $ is compensated by the rapid decrease to zero of the Fourier transform of the wavelet. 
\item The interpretation of the coherence (\ref{eq-coherence:wt})  is difficult here. Indeed, it is   complex valued, a property which is not natural for a coherence.  This comes from the fact that the quantities ${S}^{jj}_{a_1,a_2}(\omega)$ are not power spectral densities but cross-spectral densities (cross-spectral density between two different scales of the wavelet transform of one signal). Thus, to interpret correctly the coherence, we should look at one scale only, in which case we recover the coherence evaluated in the usual spectral domain. And this result is logical since the usual coherence is independent of the frequency. 
\item Setting $a_1=a_2$ and $j=k$ in the expression of the cross-spectral density, we recover the usual result of the power spectral density at one scale of the wavelet transform of a scalar fBm. The proof proposed here is a natural  extension of the proof  found in \citet{KatoM99}.

\item The derivation of the analytic form of the cross-spectral density is  easy if we use generalized functions (or Schwartz distributions). Indeed, from \citet{GelfS64}, we know that the Fourier transforms of $|v|^\alpha$ and $|v|^\alpha \sign(v)$ are respectively given by $-2 \Gamma(\alpha+1)|\omega|^{-\alpha-1} \sin (\pi\alpha/2)$
and $-2 {\bf i} \Gamma(\alpha+1)|\omega|^{-\alpha-1} \sign(\omega) \sin (\pi\alpha/2)$. Then, in the proof below, the calculation of  $T= \int_\R w_{jk}(v) \overline{\Gamma_\psi}(v) dv $ can be done using  Parseval equality. However, the theoretical background required and hidden in the calculation  is far more involved than the basics we have used in the proof (see \citet{GelfS64}). 
\end{enumerate}

\vspace{.5cm}


\begin{proof} $(i)$ We recall that under [C1] and [C2(2)], Equation~(\ref{eq-wjkgpsi}) holds, that is $E[ d_{a_1,b_1}^j  \overline{d_{a_2,b_2}^k} ] = -\frac{\sigma_j\sigma_k}2 T $ with $T:=\int_\R w_{jk}(v)\overline{\Gamma_\psi}(v)dv$. Furthermore, note that the Fourier transforms of $\psi_{a,b}$ and $\overline{\Gamma_\psi}(v)$ exist and are equal respectively to $\sqrt{a} \widehat{\psi}(a\omega)e^{-{\bf i}\omega b}$ and to 
\begin{equation}\label{def-q}
q(\omega):= \widehat{\overline{\Gamma_\psi}}(\omega) = \int_\R \overline{\Gamma_\psi}(v)e^{-i\omega v} dv =  \sqrt{a_1a_2} \overline{\widehat{\psi}(a_1\omega)} \widehat{\psi}(a_2\omega) e^{{\bf i}\omega h}.
\end{equation}
Now, let us split the proof into two cases. \\

\underline{Case 1.} $\alpha:=H_j+H_k \neq 1$.\\

When $j=k$, at this step,  \citet{KatoM99} have used the representation of $|v|^\alpha$ obtained by \citet{VonBE65}. We have obtained a similar representation for the function $\sign(v)|v|^\alpha$ for $\alpha \in (0,2)\setminus \{1\}$ (see Equations~(\ref{bahr}) and~(\ref{bahr-sign}) in Lemma~\ref{lem-bahr}). We have by Fubini's theorem and under Assumption [C3] (with $M\geq 2$).
\begin{eqnarray*}
T &=& \int_\R (\rho_{jk} -\eta_{jk} \; \sign(v))|v|^\alpha \overline{\Gamma_\psi}(v) dv \\
&=& \frac{\Gamma(\alpha+1)}\pi \int_\R |\omega|^{-\alpha-1} \int_\R  \bigg( \rho_{jk} \sin(\pi\alpha/2) (1-\cos(\omega v)) - \eta_{jk}\cos(\pi \alpha/2) \sign(\omega)(\sin(\omega v)-g_\alpha(\omega v)) \bigg)\times \\
&& \hspace*{4cm}\overline{\Gamma_\psi}(v) \; dv \; d\omega\\
&=& \frac{\Gamma(\alpha+1)}{\pi} \int_\R |\omega|^{-\alpha-1} \left(-\rho_{jk}\sin(\pi \alpha/2) \left( \frac{q(-\omega) + q(\omega)}2 \right) - \eta_{jk} \cos(\pi\alpha/2) \sign(\omega)\left(\frac{q(-\omega)-q(\omega)}{2{\bf i}} \right) \right)d\omega \\
&=&  - \frac{\Gamma(\alpha+1)}{\pi} \int_\R |\omega|^{-\alpha-1} \bigg( \rho_{jk} \sin(\pi \alpha/2)\; + {\bf i}  \;\eta_{jk} \cos(\pi\alpha/2) \sign(\omega) \bigg) q(\omega) d\omega\\
&=& - \frac{\Gamma(\alpha+1)}{\pi} \int_\R |\omega|^{-\alpha-1} \zeta_{jk}(\omega) q(\omega) d\omega.
\end{eqnarray*}
Note that the condition $M\geq 2$ is required for $\alpha>1$. For $\alpha<1$, $M\geq 1$ is a sufficient condition. These conditions allow us to show that  the contributions  $\int \overline{\Gamma_\psi}(v) dv $ and $\int_\R g_\alpha(\omega v)\overline{\Gamma_\psi}(v) dv$ are equal to zero. Now, using~(\ref{def-q}) we obtain

\begin{eqnarray*}
E[ d_{a_1,b_1}^j  \overline{d_{a_2,b_2}^k} ] &=& \sqrt{a_1a_2} \; \sigma_j\sigma_k \Gamma(\alpha+1) \times \frac1{2\pi}   \int_\R   \underbrace{|\omega|^{-\alpha-1}\zeta_{jk}(\omega) \overline{\widehat{\psi}(a_1\omega)} \widehat{\psi}(a_2\omega) }_{=:P(\omega)} e^{i\omega h} d\omega. 
\end{eqnarray*}
By using Bochner's theorem, the proof will be done, if one proves that the function $P(\cdot)$ is integrable. 
Let us prove this last assertion. Under [C2(M)], $t^k\psi(t)\in L^1$ for $k=0,\ldots,M$. Therefore, $\widehat{\psi}$ is a $M$ times continuous and differentiable function. Using a Taylor expansion
$$
\widehat{\psi}(\omega)= \sum_{k=0}^{M-1} \omega^k \widehat{\psi}^{(k)}(\omega) + \omega^M \widehat{\psi}^{(M)}(\widetilde{\omega}) = \omega^M \widehat{\psi}^{(M)}(\widetilde{\omega}), \mbox{ with } \widetilde\omega \in [0 \wedge \omega,0\vee \omega],$$
under [C2(M)]. And since $\psi^{(M)}$ is continuous at zero, $\widehat{\psi}(\omega)\sim \omega^M\widehat{\psi}^{(M)}(0)$ 
as $\omega\to 0$. Then as $\omega\to 0$:
\begin{equation}\label{eq-P}
P(\omega) \sim \zeta_{jk}(\omega)  |\omega|^{2M-1-\alpha} (a_1a_2)^M |\widehat{\psi}^{(M)}(0)|^2.
\end{equation}
As a consequence, for $M\geq 2$, $P$ is continuous at zero and $\lim_{\omega \to 0^\pm} P(\omega) =0$. Therefore for $\varepsilon>0$, $P$ is integrable on the interval $[-\varepsilon,\varepsilon]$ as a continuous function on this interval.
Finally (with $c^\vee:=|\rho_{jk}|+|\eta_{kj}|$),
\begin{eqnarray*}
\int_{|\omega|\geq \varepsilon} |P(\omega)| &\leq & c^\vee  \Big( a_1^\alpha\int_{|\omega|\geq a_1 \varepsilon } \frac{|\widehat{\psi}(\omega)|^2}{|\omega|^{\alpha+1} }d\omega\Big)^{1/2}
 \Big( a_2^\alpha\int_{|\omega|\geq a_2 \varepsilon }\frac{|\widehat{\psi}(\omega)|^2}{|\omega|^{\alpha+1}} d\omega\Big)^{1/2} \\
&\leq & \frac{c^\vee}{\varepsilon^{\alpha}} \int_\R \frac{|\widehat{\psi}(\omega)|^2}{|\omega|} d\omega <+\infty,
\end{eqnarray*}
under [C1]. Hence, $P(\cdot)\in L^1$ and Bochner's theorem may be applied.\\

\noindent \underline{Case 2.} $H_j+H_k= 1$.\\

We start with the representation of $v\log|v|$ given by (\ref{hlogh}).
\begin{eqnarray*}
w_{jk}(v)&=& \rho_{jk} |v|+ \eta_{jk} v\log|v| = \lim_{\alpha \to 1^-} \rho_{jk}|v|^\alpha +  \eta_{jk} v\log|v| \\
&=& \lim_{\alpha \to 1^-} \frac1{2\pi} \int_{\R} \frac{2\rho_{jk}(1-\cos(\omega v)) - \pi \eta_{jk} \sign(\omega) \sin(\omega v)}{|\omega|^{\alpha+1}} d\omega.
\end{eqnarray*}
Now, we derive the computation of the term $T:=\int_\R w_{jk}(v)\overline{\Gamma}_\psi(-v)dv$, similarly as the previous case. Using dominated convergence theorem and Fubini's theorem, 
\begin{eqnarray*}
T &=& \frac1{2\pi} \int_\R \left( \lim_{\alpha \to 1^-} \int_{\R} \frac{2\rho_{jk}(1-\cos(\omega v)) - \pi \eta_{jk} \sign(\omega) \sin(\omega v)}{|\omega|^{\alpha+1}} d\omega \right) \overline{\Gamma}_\psi(v) dv \\
&=& \frac1{2\pi} \lim_{\alpha \to 1^-} \int_\R \left( \int_\R \frac{2\rho_{jk}(1-\cos(\omega v)) - \pi \eta_{jk} \sign(\omega) \sin(\omega v)}{|\omega|^{\alpha+1}} \overline{\Gamma}_\psi(v) dv \right) d\omega \\
&=& \frac1{2\pi}\lim_{\alpha \to 1^-} \int_\R \left(-2\rho_{jk} \left(\frac{q(-\omega)+q(\omega)}{2}\right) - \pi \eta_{jk} \sign(\omega) \left(\frac{q(-\omega)-q(\omega)}{2{\bf i}} \right)\right) |\omega|^{-\alpha-1}
d\omega\\
&=& - \frac1{2\pi} \lim_{\alpha \to 1^-} \int_\R \frac{2\rho_{jk}+{\bf i}\pi \eta_{jk} \sign(\omega)}{|\omega|^{\alpha+1}} q(\omega) d\omega \\
&=&  - \frac1{2\pi}\lim_{\alpha \to 1^-} \int_\R |\omega|^{-\alpha-1} \left(2\rho_{jk}+{\bf i}\pi \eta_{jk} \sign(\omega)\right) \overline{\widehat{\psi}(a_1\omega)} \widehat{\psi}(a_2\omega)  e^{i\omega h}d\omega.
\end{eqnarray*}
From (\ref{eq-P}), ${|\omega|^{-\alpha-1}}\overline{\widehat{\psi}(a_1\omega)} \widehat{\psi}(a_2\omega)$ is an integrable function for all $\alpha\in (0,2)$. Therefore, the integral and the limit may be interchanged. Therefore, we obtain
$$
E[ d_{a_1,b+h}^j  \overline{d_{a_2,b}^k} ] =\sqrt{a_1 a_2}\sigma_j\sigma_k\times \frac{1}{2 \pi} 
\int_\R \frac{\rho_{jk}+{\bf i}\frac\pi2 \eta_{jk} \sign(\omega)}{|\omega|^2} \overline{\widehat{\psi}(a_1\omega)}{\widehat{\psi}(a_2\omega)}   e^{i\omega h}d\omega,
$$
and Bochner's theorem can be applied.\\
$(ii)$ is derived from (\ref{eq-P}).
\end{proof}\\

 
\section{Bahr and Essen type representations for the functions $\sign(v)|v|^\alpha$, $v_+^\alpha$ and~$v_-^\alpha$}

In 1965, von Bahr and Essen  have obtained the following representation theorem for $|v|^\alpha$ for $\alpha \in (0,2)$: 
\begin{equation}\label{bahr}
|v|^\alpha = \frac{\Gamma(\alpha+1) \sin(\pi \alpha/2)}{\pi} \int_\R \frac{1-\cos(\omega v)}{|\omega|^{\alpha+1}}d\omega.
\end{equation}
The following lemma provides a similar representation for $\sign(v)|v|^\alpha$, $v_+^\alpha=v^\alpha\mathbf{1}_{\R^+}(v)$ and $v_-^\alpha=(-v)^\alpha \mathbf{1}_{\R^-}(v)$.

\begin{lemma} \label{lem-bahr} Let  $\alpha \in (0,2)\setminus \{1\}$  and let $g_\alpha:\R \to \R$ the function which equals zero when $\alpha\in (0,1)$ and which is the identity function when $\alpha\in (1,2)$, then  we have
\begin{eqnarray}
\sign(v)|v|^\alpha &=& \frac{\Gamma(\alpha+1)\cos(\pi \alpha/2)}\pi  \int_\R \frac{\sign(\omega)\left(\sin(\omega v)-g_\alpha(\omega v) \right)}{|\omega|^{\alpha+1}} d\omega, \label{bahr-sign}\\
v_+^\alpha &=& \frac{\Gamma(\alpha+1)}{2\pi} \int_\R  \frac{\sin\left(\pi\frac{\alpha}2\right) \left( 1-\cos(\omega v)\right) + \cos\left(\pi\frac{\alpha}2\right) \sign(\omega) \left( \sin(\omega v) -g_\alpha(\omega v) \right)}{|\omega|^{\alpha+1}} d\omega, \nonumber \\
v_-^\alpha &=& \frac{\Gamma(\alpha+1)}{2\pi} \int_\R  \frac{\sin\left(\pi\frac{\alpha}2\right) \left( 1-\cos(\omega v)\right) - \cos\left(\pi\frac{\alpha}2\right)\sign(\omega) \left( \sin(\omega v) -g(\omega v) \right)}{|\omega|^{\alpha+1}} d\omega. \nonumber 
\end{eqnarray}
\end{lemma}
The representations of $v_+^\alpha$ and $v_-^\alpha$ are obtained from~(\ref{bahr}) and~(\ref{bahr-sign}) noticing that 
$$v_+^\alpha= \frac12 \left( |v|^\alpha+ \sign(v) |v|^\alpha\right) \mbox{ and } v_-^\alpha= \frac12 \left( |v|^\alpha- \sign(v) |v|^\alpha\right).$$ 

\begin{proof}
Let $\alpha \in (0,1)$, then from (\ref{bahr}) and properties of the function $\Gamma$
$$\frac1{\alpha+1} |v|^{\alpha+1} = \frac{\Gamma(\alpha+1)}\pi \cos(\pi\alpha/2) \int_\R \frac{1-\cos(\omega v)}{|\omega|^{\alpha+2}}d\omega.$$
Since $\int_\R |\omega|^{-\alpha-1} |\sin(\omega v)|<+\infty$ for $\alpha\in (0,1)$, we can differentiate this integral with respect to $v$ to obtain
\begin{equation}\label{tmp}
\sign(v)|v|^\alpha = \frac{\Gamma(\alpha+1)}\pi \cos(\pi\alpha/2) \int_\R \frac{\sign(\omega) \sin(\omega v)}{|\omega|^{\alpha+1}}d\omega.
\end{equation}
When $\alpha\in (1,2)$, then from (\ref{bahr}) and properties of the function $\Gamma$
$$\alpha |v|^{\alpha-1} = \frac{\Gamma(\alpha+1)}\pi (-\cos(\pi\alpha/2)) \int_\R \frac{1-\cos(\omega v)}{|\omega|^\alpha}d\omega. 
$$
Since $\int_\R |\omega|^{-\alpha-1}|\sin(\omega v)-\omega v|d\omega <+\infty$ for $\alpha\in (1,2)$, we can take the primitive of the last equation to get
\begin{eqnarray*}
\sign(v)|v|^\alpha &=& \frac{\Gamma(\alpha+1)}\pi \cos(\pi\alpha/2) \int_\R \frac{\sin(\omega v)/\omega-v}{|\omega|^{\alpha}}d\omega \\
&=& \frac{\Gamma(\alpha+1)}\pi \cos(\pi\alpha/2) \int_\R \frac{\sign(\omega) (\sin(\omega v)-\omega v)}{|\omega|^{\alpha+1}}d\omega, 
\end{eqnarray*}
which ends the proof.
\end{proof}

Let $\alpha \in (0,1)$, then by differentiating (\ref{tmp}) with respect to $\alpha$ and taking the limit as $\alpha\to 1^-$, we may obtain
\begin{equation}\label{hlogh}
\sign(h)|h|\log|h| = h\log|h| = \lim_{\alpha \to 1^-} - \frac12\int_\R  \frac{\sign(\omega) \sin(\omega v)}{|\omega|^{\alpha+1}}d\omega.
\end{equation}

\footnotesize
\bibliographystyle{plainnat}  
\bibliography{jtsa}

\newpage

\end{document}